\documentclass{amsart}

\usepackage{amscd,amssymb,amsthm,amsmath,amsrefs}
\usepackage[matrix,arrow]{xy}
\usepackage{mathrsfs}
\usepackage{enumerate}

\newtheorem{thm}{Theorem}[section]

\newtheorem{lemma}[thm]{Lemma}

\newtheorem{cor}[thm]{Corollary}

\theoremstyle{definition}

\newtheorem{remark}[thm]{Remark}

\numberwithin{equation}{section}

\newcommand{\Sym}{\operatorname{Sym}}

\newcommand{\Spec}{\operatorname{Spec}}

\newcommand{\E}{\mathscr{E}}

\renewcommand{\epsilon}{\varepsilon}
\renewcommand{\L}{\mathscr{L}}
\renewcommand{\O}{\mathscr{O}}

\newcommand{\op}{\operatorname}
\newcommand{\ms}{\mathscr}

\newcommand{\mb}{\mathbb}
\newcommand{\mbf}{\mathbf}
\newcommand{\mr}{\mathrm}
\newcommand{\tr}{\textrm}

\newcommand{\wt}{\widetilde}
\newcommand{\wb}{\overline}

\newcommand{\iso}{\xrightarrow{\sim}}

\newcommand{\xto}{\xrightarrow}

\newcommand{\into}{\hookrightarrow}

\newcommand{\colim}{\varinjlim}

\begin{document}

\title{$G$-torsors over a Dedekind scheme}
\author{Michael Broshi}
\address{Wellesley College, Department of Mathematics, 106 Central St, Wellesley, MA 02481}
\email{mbroshi@wellesley.edu}
\date{\today}

\begin{abstract}
We prove the equivalence of three ``points of view" on 
the notion of a $G$-torsor when the base scheme is 
a Dedekind scheme, generalizing known results when the base
is a field.  The two main tools that we generalize
are Chevalley's theorem on semi-invariants (cf. \cite{Bor}*{II.5.1})
and a Tannakian description of $G$-torsors
given by Nori and Saavedra (cf. \cite{Nori}*{Sec. 2} and \cite{Sa}*{II.4.2}). 
As an application, we show that the fibered category of $G$-torsors
on a regular proper curve over a field $k$ is an Artin
stack locally of finite presentation over $k$.
\end{abstract}

\maketitle

\section{Introduction}

Let us first fix some notation.
We fix a Dedekind scheme $X$ (the base scheme).  That is, 
$X$ is a scheme that has a finite affine open cover by the spectra of Dedekind domains.
Unless stated otherwise, any unadorned product 
is assumed to be over $X$, and for two $X$-schemes $Y$ and $T$ 
we often write $Y_{T}=Y\times T=Y\times_X T$.
If $Y$ is a scheme over $X$, we use the ``functor of points
notation" and write $y\in Y$ to denote a morphism $y:T\to Y$
of schemes over $X$.  In the same spirit, if $V$
is a locally free $\O_X$-module of finite rank, we denote also by $V$
the functor $V: T\mapsto V\otimes \O_T$, for $T$ an $X$-scheme.  This
functor is represented by $\Spec{(\Sym{V^*})}$,  where
$V^*=\op{\ms{H}\hspace{-3pt}\textit{om}}_{\O_X}{(V,\O_X)}$ denotes the dual of $V$.
For any $Y$, if $M$ is an $\O_Y$-module and $N\subset M$ is an 
$\O_Y$-submodule, we say $N$ is \emph{locally split (in $M$)} if $N$
is Zariski-locally on $Y$ a direct summand of $M$.

We fix $G$ a \emph{flat algebraic group} over $X$, 
by which we mean a flat, affine group scheme of finite type over $X$.  
By a \emph{representation of $G$}, 
we mean a finite rank, locally free $\O_X$-module $V$ with a linear
$G$-action (for details, the reader is referred to
\S\ref{coalgebras} below).  If $Y$ is an $X$-scheme, 
a \emph{$G_Y$-torsor} is a scheme $P$ faithfully flat 
and affine over $Y$, provided 
with a right $G_Y$-action such that the following two
conditions hold:
\begin{enumerate}
	\item[(i)] The map $P\to Y$ is $G_Y$-invariant.
	
	\item[(ii)]  The natural map
	\[
	  P\times_Y G_Y\to P\times_Y P; \ (p,g)\mapsto (p,pg)
	\]
	is an isomorphism.
\end{enumerate}

It follows from faithfully flat descent (\cite{EGAIV}*{2.7.1}) that
a $G_Y$-torsor is also finitely presented over $Y$, since $G$ is finitely
presented over $X$. A map $P\to P'$ of $G_Y$-torsors is a 
$G_Y$-equivariant map of $Y$-schemes.
A \emph{trivial $G_Y$-torsor} is a $G_Y$-torsor $P\to Y$ 
that is isomorphic as a $G_Y$-torsor to the projection
map $Y\times G\to Y$.  Given this terminology,
condition (ii) is equivalent to:

\begin{enumerate}[(ii$^{\prime}$)]
	\item  The map $P\to Y$ admits a section
	fppf-locally on $Y$.
\end{enumerate}

Let $X_{\mr{Zar}}$ denote the small Zariski site on $X$, that is, the 
category whose objects are open subsets $U\subset X$
and whose morphisms are inclusions.
Denote by $\op{\mbf{Rep}}{G}$ the fibered category over $X_{\mr{Zar}}$
where for an object $U$ in $X_{\mr{Zar}}$, 
$\op{\mbf{Rep}}{G}(U)=\op{Rep}_U{G}$ is the category
of representations of $G_U$ on locally free $\O_U$-modules
of finite rank.  For a scheme $Y$ over $X$, let $\mbf{Bun}_Y$
denote the fibered category over $X_{\mr{Zar}}$ 
where for an object $U$ in $X_{\mr{Zar}}$, 
$\mbf{Bun}_Y(U)=\mr{Bun}_{Y_U}$ is the category of 
all finite rank vector bundles on $Y_U$.  Both $\op{\mbf{Rep}}{G}$
and $\mbf{Bun}_Y$ are tensor categories (as described in \S\ref{categorical_viewpoint}),
and by a \emph{tensor functor} $F:\op{\mbf{Rep}}{G}\to \mbf{Bun}_Y$ we mean
a functor of fibered categories respecting the tensor structure.

Let $V$ be a representation of $G$, $\{X_1,\dots,X_r\}$ the (nonempty) connected
components of $X$ and $\mbf{i}=(i_1,\dots,i_r)$ a sequence of natural numbers.
We denote by $\bigwedge^{\mbf{i}}V$
the vector bundle such that $\bigwedge^{\mbf{i}}V|X_k =\bigwedge^{i_k}V|X_k$,
for $k=1,\dots,r$.  We denote by $t(V)$ some 
finite iteration of the operations 
$\otimes$, $\bigwedge^{\mbf{i}}$, $\mr{Sym}^{j}$, $\oplus$, and $(\cdot)^*$.
We call such an iteration a \emph{tensorial construction}.  We
remark that a tensor functor always respects the operations $\otimes$, 
$\oplus$ and $(\cdot)^*$, but need not respect $\bigwedge^{\mbf{i}}$ or $\mr{Sym}^{j}$.
However, it is a consequence of Theorem \ref{functors=torsors} that
if $Y$ is faithfully flat over $X$, and $F:\op{\mbf{Rep}}{G}\to \mbf{Bun}_Y$
is a tensor functor that is exact and faithful on the fibers over $X_{\mr{Zar}}$,
then $F$ respects all tensorial constructions.

If $V$ is a vector bundle on $X$, and $L\subset V$
is a locally split line bundle, we denote by $\op{\underline{Aut}}{(V,L)}$
the representable functor whose $T$-points are
automorphisms $f$ of $V\otimes\O_T$ such that $f(L\otimes\O_T)=L\otimes\O_T$.
We now state our main theorems.

\begin{thm}\label{maintheorem1}
	Let $G$ be a flat algebraic group over a Dedekind scheme $X$.
	There is a representation $V$ of $G$, a tensorial construction $t(V)$, and 
	a locally split line bundle $L\subset t(V)$,
	such that $G\iso\op{\underline{Aut}}{(V,L)}$.  
\end{thm}

\begin{proof}
	This is Theorem~\ref{rep-line}.
\end{proof}

\begin{thm}\label{maintheorem2}
	Let $G$ and $X$ be as above. Let $Y$ be a scheme
	faithfully flat over $X$.  There is a natural equivalence 
	that is functorial in $Y$ of the following groupoids:
	\begin{enumerate}

		\item[(i)] the groupoid of $G_Y$-torsors;
		
		\item[(ii)] the groupoid of tensor functors 
		$F:\op{\mbf{Rep}}{G}\to\mbf{Bun}_Y$ that on each
		fiber over $X_{\mr{Zar}}$ are faithful and exact.

	\end{enumerate}
\end{thm}

\begin{proof}
	This is Theorem~\ref{functors=torsors}  
	(see also Remark \ref{F(G)_explanation} for an explanation of notation.) 
\end{proof}

We can immediately state a corollary to Theorem \ref{maintheorem1},
for which we make the following defintion.  Let $V$ be a vector
bundle on $X$, $t(V)$ a tensorial construction and $L\subset t(V)$ a line bundle.
For an $X$-scheme $Y$, we define a  \emph{$Y$-twist of $(V,L)$}, 
to be a pair $(\E,\L)$ consisting of a locally free sheaf $\E$ on $Y$ 
provided with a locally split line bundle
$\L\subset t(\E)$ that is fppf-locally isomorphic as a pair to 
$(V,L)$.  That is, there is an fppf cover $Y'\to Y$ 
and an isomorphism $f:\E_{Y'}\iso V_{Y'}$
such that $f(\L_{Y'})=L_{Y'}$.  
An isomorphism of $Y$-twists $f:(\E,\L)\to(\E',\L')$  is an 
isomorphism of vector bundles $f:\E\to \E'$ such that $f(\L)=\L'$.

\begin{cor}\label{torsors=twists}
	Let $G$ and $X$ be as above. Fix a pair $(V,L)$ as in
	Theorem~\ref{maintheorem1} so that $G\iso\op{\underline{Aut}}{(V,L)}$.  
	For any scheme $Y$ over 
	$X$, there is a natural equivalence that is functorial
	in $Y$ of the following groupoids:
	
	\begin{enumerate}[(i)]
	
		\item the groupoid of $G_Y$ torsors;
		
		\item the groupoid of $Y$-twists of $(V,L)$.
		
	\end{enumerate}
	
\end{cor}

\begin{proof}
	This is a standard construction.  Given a $G_Y$-torsor $P$ and a representation
	$W$ of $G$, we can form the associated vector bundle 
	\[
		P\times^G W:=P\times W/(pg,w)\sim (p,g^{-1}w).
	\]
	Note that this construction respects tensorial constructions (see the proof
	of Lemma \ref{F_P_is_tensor} for details).
	
	Let a $G_Y$-torsor $P$ be given.  Define 
	$\E = P\times^G V$ and $\L=P\times^G L$.
	Then it is straightforward to check that $(\E,\L)$ is a $Y$-twist of $(V,L)$
	
	For a quasi-inverse, given $(\E,\L)$, we get a $G_Y$-torsor by 
	considering the associated ``frame bundle'' 
	$P = \op{\underline{Isom}}{((V_Y,L_Y),(\E,\L))}$.
\end{proof}
	
\begin{remark}
	Combining the equivalences stated in Theorem~\ref{maintheorem2} and
	Corollary \ref{torsors=twists}, we get an equivalence from the groupoid
	of functors as in Theorem~\ref{maintheorem2} and the groupoid
	of $Y$-twists of $(V,L)$.  This has a simple description.
	Namely, it is given by $F\mapsto (F(V),F(L))$.
	
	To see this, given a functor $F:\op{\mbf{Rep}}{G}\to\mbf{Bun}_Y$, 
	the equivalence in Theorem \ref{maintheorem2} assigns to $F$ the $G$-torsor
	$F(G)$ (the notation is explained in Remark \ref{F(G)_explanation}).
	Corollary \ref{torsors=twists} then assigns to $F(G)$ the pair
	$(F(G)\times^G V, F(G)\times^G L)$.  There is a map 
	$F(G)\times^G V\to F(V)$ induced by applying $F$ to the $G$-map
	$G\times V_0\to V$ (where $V_0$ is $V$ provided with the trivial $G$-action).
	That this gives a well-defined isomorphism $(F(G)\times^G V, F(G)\times^G L)
	\iso (F(V),F(L))$ is shown in the proof of Theorem~\ref{functors=torsors}.
\end{remark}

As we mentioned in the abstract, Theorems~\ref{maintheorem1} and
\ref{maintheorem2} were known when the base is a field.
Furthermore, the idea of confining oneself to locally free, finite rank
representations of $G$ (rather than all quasicoherent sheaves with $G$-action)
over Dedekind schemes is already
present in Saavedra's book on Tannakian categories \cite{Sa}.  
Nonetheless, the equivalence in Theorem \ref{maintheorem2} is
only proven there when the base is a field (cf. \cite{Sa}*{II.4.2.2}).

Finally, we remark that the formalism involving fibered categories over 
the Zariski site on $X$ used in Theorem \ref{maintheorem2}
is not necessary when $X$ is affine.  In that case, one need
only consider exact, faithful tensor functors $F:\op{Rep}{G}\to\mr{Bun}_X$.

\textbf{Acknowledgements}: I would like to thank Madhav Nori for
many helpful discussions, and Torsten Wedhorn and Philipp Gross for
their useful correspondence.  
This paper stems from a result from my PhD thesis, and I am
deeply grateful to my advisor, Mark Kisin, for his help and support along the way.
I am indebted to Brian Conrad who generously read earlier drafts, and
gave numerous comments and suggestions.
Finally, I am happy to thank the referee who gave a very
careful reading that caught several infelicities and
improved the exposition.

\section{Application to the moduli of $G$-torsors}\label{moduli}

Before proceeding with the proof of Theorem~\ref{maintheorem1}, 
we give an application to the stack of $G$-torsors over a curve.
By an \emph{Artin stack}, we mean an algebraic stack
as defined in \cite{LMB}*{4.1}.  In particular, we
assume that an Artin stack has a separated and quasicompact
diagonal.
For this section only, let $k$ be a field, and assume that $X$ is a 
connected, regular, proper curve over $k$.  
In particular, $X$ is a Dedekind scheme.
We also assume for this section that $G$ \emph{has connected generic fibre}.
Finally, for this section only we use the convention that for $k$-schemes 
$Y$ and $T$, $Y_T=Y\times_{\Spec{k}}T$. 

Let $G\mr{Tor}_X$ denote
the fibered category that assigns to a $k$-scheme $T$
the groupoid of $G_{X_T}$-torsors.  The goal of this section is to
prove the following theorem.  We are grateful to Brian Conrad
for pointing out this application of Theorem~\ref{maintheorem1}.

\begin{thm}\label{Artin}
The fibered category $G\mr{Tor}_{X}$ is an Artin stack,
locally of finite presentation over $k$.
\end{thm}

We recall the following definition from \cite{R-G}*{3.3.3},
a key input into the proof of the theorem,
although the reader can take the statements of the subsequent theorem  and 
lemmas as a black box.  Let $S$ be a scheme and
$T$ a scheme locally of finite presentation over $S$.  We 
define the \emph{relative associated primes} of $T$ over $S$, denoted
$\op{Ass}{(T/S)}$, by 
\[
	\op{Ass}{(T/S)}=\bigcup_{s\in S}\op{Ass}{(T_s)}.
\]
For a point $s\in S$, 
denote by $(\wt{S},\wt{s})$ the henselization of the
pair $(S,s)$, and let $\wt{T}=T\times_S\wt{S}$.  We say that $T$
is \emph{pure along $T_s$} if 
for each element $\wt{t}\in\op{Ass}{(\wt{T}/\wt{S})}$, 
the closure
of $\wt{t}$ in $\wt{T}$ meets $\wt{T}_{\wt{s}}$.
We say that $T$ is $S$-\emph{pure} (or that the map $T\to S$
is \emph{pure}) if it is pure along 
$T_s$ for each $s\in S$.  

A simple example of a map that is not pure
is given by $S=\Spec{R}$ for $R$ a complete DVR, $T=\Spec{K}$ 
where $K$ is the fraction field of $R$ and $T\to S$ the natural
inclusion.  Then $T_s$ is in fact empty for $s$ the closed point of $S$.

The reason why we introduce this notion of purity is that pure
maps have ``flattening stratifications." More precisely, 
we have the following theorem.

\begin{thm}\label{pure}
Suppose that $T\to S$ is pure.  Then there is a
monomorphism $Z\into S$ that is locally of finite presentation
such that for any $S$-scheme $S'$,  $T\times_S S'\to S'$
is flat if and only if $S'\to S$ factors through $Z$.
\end{thm}

\begin{proof}
This is \cite{R-G}*{I.4.3.1}.
\end{proof}

\begin{lemma}\label{G-pure}
With $G$ and $X$ as above, $G$ is $X$-pure.
\end{lemma}

\begin{proof}
Let $\xi\in X$ be the generic point of $X$.  By assumption $G_{\xi}$ is 
connected, so it is in fact geometrically irreducible 
by \cite{SGA3}*{VI${}_{\mr{A}}$ 2.4}.  
By \cite{EGAIV}*{2.3.7}, since $G$ is
flat over $X$, and $X$ is irreducible, 
the image  of $G_{\xi}$ in $G$ is dense.  In particular, since $G_{\xi}$ is irreducible so 
is $G$.  Since $X$ has
$\xi$ as its unique associated prime, 
$\op{Ass}{G}=\op{Ass}{G_{\xi}}$ by the $X$-flatness of $G$
(see \cite{EGAIV}*{3.3.1}, which describes associated primes along fibers).
Let $\eta\in G$ be its generic point.  
We claim that $\op{Ass}{G_{\xi}}=\{\eta\}$.  
Suppose on the contrary that $Z\subset G_{\xi}$ is an embedded component.
In particular $\op{dim}{G_{\xi}}>0$.
Denote by $\wb{\xi}$ an algebraic closure of $\xi$.  Then
$Z_{\wb{\xi}}\subset G_{\wb{\xi}}$ is a union of finitely
many embedded components.
Furthermore, for each closed point $g\in G_{\wb{\xi}}$, $gZ_{\wb{\xi}}$
is also a union of finitely many embedded components.  Since $G_{\wb{\xi}}$
is irreducible, there must be infinitely many distinct closed sets
amongst the pairwise disjoint $\{gZ_{\wb{\xi}}\}_{g\in G}$.  
But, this is a contradiction since
$G_{\wb{\xi}}$ is of finite type over $\wb{\xi}$ hence has
only finitely many associated primes.

Thus far, we have concluded that $G$ is an irreducible scheme over
$X$, and its generic point $\eta\in G$ is its unique associated prime.
To show $G$ is pure over a closed point $x\in X$ we may replace 
$X$ by $\Spec{\O_{X,x}}$. So, we may assume that $X$ is the spectrum
of a DVR with closed point $x$ ($G$ is still irreducible and its
generic point is its unique associated prime after this base change).  
Let $(\wt{X},\wt{x})$ be the henselization of $(X,x)$.
Then $\wt{X}$
has its generic point at its unique associated prime.
It then follows as above that $\wt{G}:=G\times_X \wt{X}$ also has its
generic point as its unique associated prime.  Thus,
$\op{Ass}{(\wt{G}/\wt{X})}$ consists of the generic point of $\wt{G}$
together with points on $\wt{G}_{\wt{x}}$ (in fact just the generic points of
the latter, but this is not needed).  In particular, the closures of these
points in $\wt{G}$ meet $\wt{G}_{\wt{x}}$.  This shows
that $G$ is pure along $G_x$ for each closed point $x\in X$,
and it is straightforward to check that $G$ is pure along $G_{\xi}$ as well.
Hence, $G$ is pure over $X$, as claimed.
\end{proof}
	
\begin{lemma}\label{pure_pullback}
Let $T\to S$ be locally of finite presentation.
If $S'\to S$ is fppf, then $T\times_S S'\to S'$ is 
flat and pure if and only if $T\to S$ is flat and pure.
\end{lemma}

\begin{proof}
For purity this is \cite{R-G}*{I.3.3.7}, and
for flatness this is \cite{EGAIV}*{2.5.1}.
\end{proof} 

\begin{lemma}\label{represent_pure_maps}
Let $\ms{I}$ and $\ms{Q}$ be an Artin stacks
over $k$, and let $f:\ms{I}\to X_{\ms{Q}}$
be representable in schemes and locally of finite presentation.  The condition
on $\ms{Q}$-schemes $T$ that $\ms{I}\times_{\ms{Q}} T\to X_T$
is flat and pure is representable by an Artin stack
locally of finite presentation over $\ms{Q}$.
\end{lemma}

\begin{proof}
Let $\ms{Z}$ denote the fibered category over $\ms{Q}$ 
where $\ms{Z}(T)\subset\ms{Q}(T)$ is the full subcategory
consisting of those objects of $\ms{Q}(T)$ for which 
$\ms{I}\times_{\ms{Q}} T\to X_T$ is flat and pure.
Using Lemma \ref{pure_pullback}, it is straightforward to 
verify that $\ms{Z}$ is a stack.  We must show
that the map $\ms{Z}\to\ms{Q}$ is representable and
locally of finite presentation.

Let $Q\to\ms{Q}$ be a smooth scheme cover, and 
let $I=\ms{I}\times_{\ms{Q}}Q$, a smooth scheme cover of $\ms{I}$.
It suffices to show that $Z=\ms{Z}\times_{\ms{Q}}Q$ is
an algebraic space, locally of finite presentation over $Q$. By definition, 
for any $k$-scheme $T$, a map $T\to Q$ lies in $Z(T)\subset Q(T)$
if and only if $I\times_Q T\to X_T$ is flat and pure.
Thus, we must represent that condition on $Q$-schemes.
We first represent the purity condition.
By \cite{R-G}*{3.3.8}, purity is an open condition,
so there is an open immersion $U'\into X_{Q}$ such that
$X_T\to X_{Q}$ factors through $U'$ if and only if 
$I\times_Q T=I\times_{X_Q}X_T$ is
pure over $X_T$.  To get an open subspace of $Q$
representing the purity condition, we take
the (closed) image of the closed complement of $U'$ under $X_{Q}\to Q$
and let $U$ be complement of that image.  It then follows that $T\to Q$
factors through $U$ if and only if $I\times_Q T$ is pure over $X_T$.

Thus, replacing $Q$ by $U$ and $I$ by the inverse 
image of $X_U$, we may assume that $I\to X_{Q}$
is pure.  In this case, by Theorem \ref{pure}, there is a representable monomorphism
$Z'\to X_{Q}$ such that $Y\to X_{Q}$ factors through $Z'$  
if and only if $I\times_{X_{Q}} Y\to Y$ is flat.  We now
want to represent the condition on $Q$-schemes $T$ that  
$X_T\to X_{Q}$ factors through $Z'$.  These are exactly
the $T$-points of the restriction of scalars $\op{Res}^{X_{Q}}_{Q}(Z')$.
By \cite{Olsson}*{1.5}, since $X_{Q}\to Q$ is a proper, flat, and locally finitely
presented, and $Z'\to X_{Q}$ is separated and locally of finite
presentation, $\op{Res}^{X_{Q}}_{Q}(Z')$ 
is an algebraic space, locally of finite presentation over $Q$.
\end{proof}

\begin{proof}[Proof of Theorem \ref{Artin}]
By Theorem \ref{maintheorem1}, we can find a representation
of $G$ on a finite rank vector bundle $V$, a tensorial
construction $t(V)$ and a locally split line bundle $L\subset t(V)$ such that
$G\iso\op{\underline{Aut}}{(V,L)}$.  We now fix such a pair $(V,L)$.
Since $X$ is connected, $V$ has constant rank $n$ for some $n\in\mb{N}$.
For any $X$-scheme $Y$, the
identification $G\iso\op{\underline{Aut}}{(V,L)}$
pulls back to $G_Y\iso\op{\underline{Aut}}{(V_Y,L_Y)}$.
Let $\mr{Bun}^n_X$ denote the stack of rank $n$ vector bundles
over $X$ (where $n$ is the rank of $V$).  
That is, to each $k$-scheme $T$, $\mr{Bun}^n_X(T)$ is
the groupoid of rank $n$ vector bundles over $X_T=X\times_{k} T$.  
By \cite{LMB}*{4.6.2.1},
$\mr{Bun}^n_X$ is an Artin stack, locally of finite presentation over $k$.  

Let $\E^{\mr{univ}}$ denote the universal rank $n$ vector
bundle on $X\times \mr{Bun}^n_X$.  Let $\ms{Q}$ denote the
relative quot scheme over $\mr{Bun}^n_X$ classifying all
rank one, locally split subbundles of $t(\E^{\mr{univ}})$ (where $t$ is the same
tensorial construction as that defining $G$).  That is, for a scheme
$T$ over $\mr{Bun}^n_X$, $\ms{Q}(T)$
is the groupoid of locally split line bundles 
$\L_{X_T}\subset t(\E^{\mr{univ}})_{X_T}=t(\E^{\mr{univ}}_{X_T})$ on $X_T$.
Since $X$ is projective over $k$, it follows from 
\cite{FGA}*{no. 221 Theorem 3.1} that $\ms{Q}\to \mr{Bun}^n_X$ is 
representable and locally of finite presentation.
Let $\L^{\mr{univ}}\subset t(\E^{\mr{univ}}_{X_{\ms{Q}}})$ denote
the universal line bundle on $X_{\ms{Q}}$.  Finally, 
let $\ms{I}$ over $X_{\ms{Q}}$ denote the fibered category, 
where for an $X_{\ms{Q}}$-scheme $T$, 
$\ms{I}(T)=\op{Isom}{((V_T,L_T),(\E^{\mr{univ}}_T,\L^{\mr{univ}}_T))}$.
Then $\ms{I}\to X_{\ms{Q}}$ is representable in schemes, 
affine and of finite presentation.

By Lemma~\ref{represent_pure_maps}, there is an Artin stack $\ms{Z}$ 
locally of finite presentation over $\ms{Q}$
representing the condition on $\ms{Q}$-schemes $T$
that $\ms{I}\times_{\ms{Q}}T$ is flat and pure over $X_T$.  In particular, 
$\ms{I}\times_{\ms{Q}}\ms{Z}$ is flat over $X_{\ms{Z}}$.  Let 
$\ms{U}'\subset X_{\ms{Q}}$
denote its open image.  Let $\ms{U}\subset\ms{Q}$ denote the
complement of the closed image
of the complement of $\ms{U}'$ under the projection
$X_{\ms{Q}}\to\ms{Q}$.  Thus, $\ms{U}$ represents the condition on $\ms{Q}$-schemes
$T$ that $\ms{I}\times_{\ms{Q}}T$ 
is flat, surjective (hence fppf since $\ms{I}\to X_{\ms{Q}}$
is finitely presented) and pure over $X_T$.  Furthermore,
we still have that $\ms{U}$ is locally of finite presentation over $\ms{Q}$.
We now show that $\ms{U}$ is naturally isomorphic to $G\mr{Tor}_X$.
By Corollary \ref{torsors=twists}, $G\mr{Tor}_X$ 
is isomorphic to the fibered category
that assigns to a $k$-scheme $T$ the groupoid of
$X_T$-twists of $(V,L)$.  It suffices to show that
$\ms{U}$ is naturally isomorphic to this latter fibered category.

Let $T$ be a $\ms{Q}$-scheme and denote by $f:X_T\to X_{\ms{Q}}$
the corresponding map.  For ease, we write 
$f^*\ms{I}$ for the pullback of $\ms{I}$ along $f$.
The map $f:X_T\to X_{\ms{Q}}$ gives rise to a pair
$(f^*\E^{\mr{univ}}_{X_{\ms{Q}}},f^*\L^{\mr{univ}})$.  
We claim that 
$(f^*\E^{\mr{univ}}_{X_{\ms{Q}}},f^*\L^{\mr{univ}})$ is an
$X_T$-twist of $(V,L)$ if and only if $T$ factors through $\ms{U}$.
First assume that $T\to\ms{Q}$ factors through $\ms{U}$.  
In particular, $f^*\ms{I}\to X_T$ is fppf.
Note that the canonical projection $f^*\ms{I}\to \ms{I}$
gives an isomorphism $(\E^{\mr{univ}}_{f^*\ms{I}},\L^{\mr{univ}}_{f^*\ms{I}})\cong
(V_{f^*\ms{I}},L_{f^*\ms{I}})$.  Thus, $f^*\ms{I}\to X_T$
gives the desired fppf cover.  Conversely, if 
$(f^*\E^{\mr{univ}}_{X_{\ms{Q}}},f^*\L^{\mr{univ}})$ is an 
$X_T$-twist of $(V,L)$, then $f^*\ms{I}$ is a $G_{X_T}$-torsor 
(cf. the proof of Corollary~\ref{torsors=twists}), and so fppf
over $X_T$. Furthermore, since $G$ is $X$-pure by
Lemma~\ref{G-pure}, it follows by Lemma~\ref{pure_pullback}
the $G_{X_T}$-torsor $f^*\ms{I}$ is $X_T$-pure.  Thus, $T$
factors through $\ms{U}$.
We conclude that $\ms{U}$ is naturally isomorphic 
to the desired fibered category, which completes
the proof.
\end{proof}

\section{Algebraic groups over Dedekind schemes}\label{coalgebras}

With notation as in the introduction, let $G$ be a flat algebraic 
group scheme over $X$.  This means that $G$ is a flat affine group scheme of finite
type over $X$.
%
Let $f:G\to X$ denote the structure map.
We will abuse notation and denote the $\O_X$-bialgebra $f_*(\O_G)$ simply by
$\O_G$.  Let $\Delta:\O_G\to \O_G\otimes \O_G$ 
denote the comultiplication map and $\epsilon:\O_G\to \O_X$ the counit.  
%
%
%
As above, if $W\subset V$ is Zariski-locally on $X$ a direct summand as
an $\O_X$-module, we will call the inclusion \emph{locally split}.
If $W$ and $V$ are (compatibly) $\O_G$-comodules, 
that the inclusion $W\subset V$ is locally split 
does not imply in general that $W\subset V$ is locally a direct summand
as an $\O_G$-comodule. 
Finally, recall that $GL(V)$ is an algebraic group scheme that is 
represented by $\Spec{(\Sym{(V\otimes V^*)}[1/\mr{det}])}$.
Our presentation of this section follows \cite{Wat}*{Chap. 3} 
and \cite{Bor}*{Chap. 5}, generalized to our current situation.

\begin{lemma}\label{directed_union}
Let $V$ be an $X$-flat quasicoherent $\O_G$-comodule.  Then 
$V$ is the direct limit of $\O_G$-comodules that are 
locally free $\O_X$-modules of finite rank.
\end{lemma}

\begin{proof}
For $X$ affine, this is the Corollary to Proposition 1.2 in \cite{Se68}.
We quickly sketch the proof in the general case as the details are the
same as in \emph{ibid}.  Since $X$ is noetherian, by \cite{EGAI}*{9.4.9}
any quasicoherent sheaf is the direct limit of its coherent subsheaves.
Since a coherent $\O_X$-submodule of $V$ is locally free, 
it suffices to show that for any coherent
submodule $W\subset V$, $W$ is contained in a coherent 
$\O_G$-subcomodule of $V$.  Let $\rho:V\to V\otimes\O_G$
denote the comodule map.  Since $\rho(W)$ is coherent, 
there is a coherent submodule $W'\subset V$
such that $\rho(W)\subset W'\otimes \O_G$.  Define
a quasicoherent $\O_X$-module $E=\rho^{-1}(W'\otimes \O_G)$.  
By working over open affines in $X$, one can show that $E\subset W'$, so
it is coherent, and $E$ is an $\O_G$-comodule (cf. \cite{Se68}*{Section 1.5}).
\end{proof}

\begin{lemma}\label{faithfulrep}
There is a representation $V$ of $G$ such that the map
$G\to GL(V)$ is a closed embedding.
\end{lemma}

\begin{proof}
Consider the regular representation $\Delta: \O_G\to \O_G\otimes \O_G$.  
By Lemma \ref{directed_union}, there 
is a locally free, finite rank $\O_G$-subcomodule $V\subset \O_G$
that locally contains a finite system of $\O_X$-algebra generators of $\O_G$.
By restricting $\Delta$ to $V$,  we have an $\O_G$-comodule 
$\rho: V\to V\otimes \O_G$.  To check the corresponding map $G\to GL(V)$ is a
closed embedding, we may assume that $X=\Spec{R}$, where $R$ is a DVR.
In this case, $V\cong R^n$, and $\O_{GL(V)}\cong R[x_{11},\dots,x_{nn}][1/\det]$.
The verification that $\O_{GL(V)}\to\O_G$ is surjective is then identical to the proof
in \cite{Wat}*{3.4}.

Namely, if we choose a basis $\{v_1,\dots,v_n\}$ of $V$ and write
$\rho(v_i) = \sum v_j \otimes a_{ij}$, then the map $\O_{GL(V)}\to\O_G$
is given by $x_{ij}\mapsto a_{ij}$.  Since 
$v_j=(\epsilon\otimes 1)\Delta(v_j)=\sum\epsilon(v_i)a_{ij}$,
the image of the map $\O_{GL(V)}\to\O_G$ contains $V$, hence 
is surjective since $V$ contains the algebra generators of $\O_G$.
\end{proof}

Let $\{X_1,\dots,X_r\}$ denote the set of (nonempty) connected
components of $X$.  Let $\ms{K}_{X_i}$ denote the stalk of $\O_{X_i}$
at the generic point of $X_i$, and write $\ms{K}_X=\prod\ms{K}_{X_i}$.
If $M$ is a locally free of finite rank $\O_X$-module, 
and $N'\subset M$ is a coherent submodule,
we call $N=(N'\otimes\ms{K}_X)\cap M\subset M\otimes\ms{K}_X$
the \emph{saturation of $N'$ in $M$}.  (The point is that $N'$
may not be a subbundle of $M$.)

\begin{lemma}\label{saturation}
Let $W$ be a representation of $G$, $U'\subset W$
a subrepresentation, and let $U$ denote the saturation of 
$U'$ in $W$. Then, $U$ is a subrepresentation of $W$ that is locally
split as an $\O_X$-module.
\end{lemma}

\begin{proof}
Since $X$ is Dedekind, it
is straightforward to check that $U$ is locally split in $W$
(say, by looking at stalks and using the elementary divisors theorem). It
remains to show that $U$ is $G$-stable.
Let $\rho: W\to W\otimes\O_G$ denote the comodule
map.  We wish to show that $\rho(U)\subset U\otimes\O_G$.
This can be checked Zariski-locally on $X$, so can 
assume that $X=\Spec{A}$ is a Dedekind domain, and 
$U/W$ is free.  To show that the image of
$U$ in $W\otimes\O_G$ is contained in $U\otimes\O_G$, we must
show the image of any element in $W\otimes\O_G$  goes to zero in $(U/W)\otimes\O_G$.
Since this latter $A$-module is flat, we can check that the image is zero
on the generic point of $\Spec{A}$.  But, over the generic point $U=U'$, so
the result follows from the $G$-stability of $U'$.
\end{proof}

\begin{lemma}\label{linalg}
Let $W$ be a finite rank vector bundle on $X$, and suppose $U\subseteq W$ 
is a locally split subbundle. Let $\mbf{d}=(d_1,\dots,d_r)$
be the sequence of ranks of $U$ on each nonempty connected component of $X$.
Define $L=\bigwedge^{\mbf{d}}U\subset\bigwedge^{\mbf{d}} W$.  Let $g\in GL(W)$.  Then 
\[
	gL=L \iff gU=U.
\]
\end{lemma}

\begin{proof}
The statement is local on $X$, so we suppose that $X=\Spec{A}$
for a Dedekind domain $A$, and that $U\subset W$ is a rank $d$ direct summand.
The direction $\Longleftarrow$ is immediate by functoriality, so we assume now that 
$gL=L$. First, note that for any $A$-algebra $B$, 
\[
	U\otimes B=\{\omega\in W\otimes B\mid \omega\wedge (L\otimes B)=0\}.
\]
If $g\in GL(W\otimes B)$ and $u\in U\otimes B$, then
\[
	gu\wedge (L\otimes B)=g(u\wedge g^{-1}(L\otimes B))=g(u\wedge L\otimes B)=0.
\]
It follows from the previous remark that $gu\in U\otimes B$, as desired.
\end{proof}

\begin{thm}\label{rep-line}
There is a representation $V$ of $G$, a tensorial
construction $t(V)$, and
a locally split line bundle $L\subset t(V)$ such that 
\[ 
	G=\{g\in GL(V)\mid gL=L\}.
\]
\end{thm}

\begin{proof}
By Lemma \ref{faithfulrep}, we can fix a representation $V$ of
$G$ such that $G\to GL(V)$ is a closed embedding.  We must now construct
$t(V)$ and $L\subset t(V)$.  We can write 
\begin{equation}\label{GL(V)}
	\ms{O}_{GL(V)}=\colim_i
	{\left(\bigoplus_{m\geq 0}\op{Sym}^m{(V\otimes V^*)}\cdot\det{}^{-i}\right)}.
\end{equation}
Identifying $G$ as a closed subgroup of $GL(V)$, 
$G$ is defined by a coherent sheaf of ideals $\ms{I}\subset \ms{O}_{GL(V)}$. 
Note that since $G$ is flat over $X$, $\ms{I}$ is saturated in $\ms{O}_{GL(V)}$.
Choose a finite open affine cover $\{X_i\}$ of $X$.
On each $X_i$, $\ms{I}|X_i$ is finitely generated in $\ms{O}_{GL(V)}|X_i$ 
as an $\O_{X_i}$-module.  Hence,  by taking integers $M$ and $N$
sufficiently large, we can ensure that
the module generators of $\ms{I}$ on each $X_i$ are contained in
\[
	t'(V)=\bigoplus_{m=0}^M\op{Sym}^m{(V\otimes V^*)}\cdot\det{}^{-N}.
\]
Let $U'=\ms{I}\cap t'(V)$. 
Let $G'=\{g\in GL(V)\mid gU'=U'\}$.  We claim that $G=G'$.
First, note that
\[
	G=\{g\in GL(V)\mid g\ms{I}=\ms{I}\}.
\]
In particular, $G\subseteq G'$.  On the other hand, 
if $g\in G'(B)$, then by definition the
induced map $(1\otimes g)\circ\Delta: U'\to\ms{O}_{GL(V)}\otimes B$
factors through $U'\otimes B$.  However, since $(1\otimes g)\circ\Delta$
is an $\O_X$-algebra map, it follows that $\ms{I}\to \ms{O}_{GL(V)}\otimes B$ factors
through $\ms{I}\otimes B$.  That is, $G'\subseteq G$, thus $G=G'$.

Let $U$ be the saturation of $U'$ in $t'(V)$. 
By Lemma~\ref{saturation}, $U$ is $G$-stable and
locally split in $t'(V)$.  It follows that $G\subseteq\{g\in GL(V)\mid gU=U\}$.
Conversely, to check that $\{g\in GL(V)\mid gU=U\}\subseteq G$,
it suffices to check on an affine cover of $X$.  Then one can
see that $\{g\in GL(V)\mid gU=U\}\subseteq G$ exactly as in the
proof of Lemma~\ref{saturation}.  Thus, $G=\{g\in GL(V)\mid gU=U\}$.
Let $\mbf{d}=(d_1,\dots,d_n)$
be the sequence of ranks of $U$ on each nonempty connected component of $X$.
Define $t(V)=\bigwedge^{\mbf{d}} t'(V)$ and 
$L=\bigwedge^{\mbf{d}}U\subset t(V)$. By Lemma~\ref{linalg},
we have that $G=\{g\in GL(V)\mid gL=L\}$, as claimed.
\end{proof}

\section{Tannakian viewpoint}\label{categorical_viewpoint}

We recall the notation from the introduction.
As usual, $G$ denotes a flat algebraic group over a Dedekind scheme $X$.
In this section, we fix a faithfully flat  $X$-scheme $Y$.   
Recall that unadorned products are fiber products over $X$ and for an 
$X$-scheme $T$, $Y_T=Y\times T=Y\times_X T$.  
For each open subscheme $U\subset X$, let $\O_U=\O_X|U$. We write 
$\op{Rep}_U{G}$ for the category of
representations of $G_U$ on finite rank, locally free $\O_U$-modules.  
Then, $\op{Rep}_U{G}$ is an $\O_U$-linear, rigid tensor category.
Here, rigid means that $\op{Rep}_{U}{G}$ has internal homs.
Of course, unless $\O_U$ is a field,
this will not be an abelian category.  Denote by $X_{\mr{Zar}}$
the small Zariski site on $X$.
Denote by $\op{\mbf{Rep}}{G}$ the fibered over $X_{\mr{Zar}}$
where for an object $U$ in $X_{\mr{Zar}}$, 
$\op{\mbf{Rep}}{G}(U)=\op{Rep}_U{G}$.  
Then $\op{\mbf{Rep}}{G}$ is a \emph{(fibered) tensor category} in the following sense: 
\begin{enumerate}[(i)]
	\item There is a monoidal structure
	\[
		\op{\mbf{Rep}}{G}\times_{X_{\mr{Zar}}}\op{\mbf{Rep}}{G}\to \op{\mbf{Rep}}{G}
	\]
	(along with associativity and commutativity constraints)
	that over each $U$ in $X_{\mr{Zar}}$ induces the usual tensor structure
	on $\op{Rep}_U{G}$.
		
	\item There is an object $1_X\in\op{Rep}_X{G}$ that pulls back to the unit
	object in $\op{Rep}_U{G}$ for each $U$ in $X_{\mr{Zar}}$.
	
	\item For each $U'\subset U$, the pullback
	map $\op{Rep}_U{G}\to\op{Rep}_{U'}{G}$ is a tensor functor.
\end{enumerate}
Let $\mbf{Bun}_Y$ denote the fibered category over $X_{\mr{Zar}}$ 
where for an object $U$ in $X_{\mr{Zar}}$, 
$\mbf{Bun}_Y(U)=\mr{Bun}_{Y_U}$ is the category of all
finite rank vector bundles on $Y_U$ (not to
be confused with $\mr{Bun}^n_Y$ in $\S$\ref{moduli}).
Then $\mbf{Bun}_Y$ is a tensor category each of whose fibers over $X_{\mr{Zar}}$
is $\O_{Y_U}$-linear and rigid.  
By a \emph{(fibered) tensor functor} 
$F:\op{\mbf{Rep}}{G}\to \mbf{Bun}_Y$, we mean a functor of fibered categories
over $X_{\mr{Zar}}$ that induces a tensor functor (in the usual sense) on each
fiber.  In particular, $F$ must respect unit objects on each fiber.

Let $P\to Y$ be a $G_Y$-torsor.  Then, for each object $U$ 
in $X_{\mr{Zar}}$, $P_U$ is a $G_{Y_U}$-torsor.
We define a functor $F_P:\op{\mbf{Rep}}{G}\to\mbf{Bun}_Y$ as
follows. For an object $U$ in $X_{\mr{Zar}}$, and $V$ in 
$\op{Rep}_U{G}$, 
\[
	F_P: V\mapsto P_{U}\times^{G_{Y_U}} (V\times_U{Y_U}) = 
	P_U\times (V\times_U {Y_U})/((p,v)\sim (pg,g^{-1}v)).
\]
Concretely, we are pushing out $P$ along the map $G\to GL(V)$ to
associate to the $G_{Y_U}$-torsor $P$ a $GL(V_U)$-torsor, that is, a
vector bundle on $Y_U$.
It is clear $F_P$ respects pullback maps, so it is a
functor of fibered categories.
When no confusion will arise, we will write $F_P(V)=P\times^G V$
for notational ease.

\begin{lemma}\label{F_P_is_tensor}
The functor $F_P$ is a tensor functor that on each fiber over $X_{\mr{Zar}}$ is 
faithful and exact.
\end{lemma}

\begin{proof}
It is clear that $F_P$ is a functor of fibered categories
over $X_{\mr{Zar}}$, so we must show it is an exact, faithful
tensor functor on each fiber. Fix an object $U$ in $X_{\mr{Zar}}$.
Since $G$ acts transitively on $P$, it is straightforward to check
from the definition that $F_P(\O_U)=\O_{Y_U}$, where $\O_U$ has the trivial 
$G_U$-action.  Thus, $F_P$ respects unit objects.
For $V$ and $W$ in $\op{Rep}_U{G}$, there is a natural map
\begin{equation}\label{associated_bundle_tensor_map}
	P\times^{G} (V\otimes W)\to (P\times^{G}V)\otimes (P\times^{G} W);
	\ (p,v\otimes w)\mapsto (p,v)\otimes (p,w),
\end{equation}
which is straightforward to check is well defined.
It suffices to check that \eqref{associated_bundle_tensor_map} is an isomorphism 
fppf-locally on $U$, so we may assume that 
$P_{U}=G_{Y_U}\times_U Y_U$ is the trivial $G_{Y_U}$-torsor.
Under the identification, $P\times^G V=V_{Y_U}$, 
\eqref{associated_bundle_tensor_map} becomes the identity map.  Hence, $F_P$
is a tensor functor.

To show that $F_P$ is exact, we must show that if 
$0\to V'\to V\to V''\to 0$ is exact, then so is 
$0\to P\times^G V'\to P\times^G V\to P\times^G V''\to 0$.
Again, we can check that this sequence is exact fppf-locally
on $Y_U$, so we can assume that $P_U=G_{Y_U}\times_U Y_U$.  We can then
identify the latter exact sequence with $0\to V'_{Y_U}\to V_{Y_U}\to V''_{Y_U}\to 0$,
which is exact since $Y$ is flat over $X$.  Finally, to show
that $F_P$ is faithful, we assume that $F_P(V)=0$.  Passing to
an fppf-cover of $Y_U$, this implies that $V_{Y_U}=0$.  Hence $V=0$
since $Y$ is faithfully flat over $X$.  
\end{proof}

\begin{remark}\label{generalize_F_P_is_tensor}
The proof above that $F_P$ respects tensor products
generalizes easily to show that in fact $F_P$ respects any
tensorial construction.
\end{remark}

Thus, $F_P$ is a tensor functor that on each fiber over $X_{\mr{Zar}}$ is faithful 
and exact. We now prove that the converse is true.
Let $F:\op{\mbf{Rep}}{G}\to\mbf{Bun}_Y$ be a tensor functor that on 
each fiber over $X_{\mr{Zar}}$ is faithful 
and exact. We show that there is a natural equivalence 
$F\iso F_P$ for a uniquely defined $G_Y$-torsor $P$.  
We closely follow the elegant presentation in \cite{Nori}*{Sec. 2},
generalizing to our current situation.  The main idea
to define $P$ is to apply $F$ to the regular representation of $G$.
Of course, this is not a finite rank representation, so we must
first suitably extend $F$.

We denote by $\op{\mbf{Rep}}^{\prime}{G}$ the 
fibered category over $X_{\mr{Zar}}$, where for each $U$ in $X_{\mr{Zar}}$,
$\op{\mbf{Rep}}^{\prime}{G}(U)=\op{Rep}^{\prime}_U{G}$ is the category of 
flat quasicoherent $\O_U$-modules
that are also $\O_{G_U}$-comodules.  Denote by $\mbf{QCoh}_Y$ the fibered
category over $X_{\mr{Zar}}$ where for each $U$ in $X_{\mr{Zar}}$, 
$\mbf{QCoh}_Y(U)=\mr{QCoh}_{Y_U}$
is the category of quasicoherent $\ms{O}_{Y_U}$-modules.

Since we will be working over open subschemes of $X$, we
will need the following slight generalization of Lemma \ref{directed_union}.

\begin{lemma}\label{generalized_directed_union}
Let $U\subset X$ be an open subscheme and let $V$ be an object
of $\op{Rep}^{\prime}_U{G}$.  Then, $V$ is the direct limit of 
its subobjects in $\op{Rep}_U{G}$. 
\end{lemma}

\begin{proof}
The proof is identical to that of Lemma \ref{directed_union}.  
One need only note that
it is still the case that any coherent $\O_U$-submodule of $V$ is 
locally free.
\end{proof}

\begin{lemma}\label{functor_extension}
The functor $F$ extends uniquely to a tensor functor 
$F:\op{\mbf{Rep}}^{\prime}{G}\to\mbf{QCoh}_Y$ such that:
\begin{enumerate}[(i)]
	\item On each fiber over $X_{\mr{Zar}}$, $F$ is exact and faithful.
	\item The extended $F$ respects direct limits.
	\item The $\O_{Y}$-module $F(\O_G)$ is faithfully flat.
\end{enumerate}
\end{lemma}

\begin{proof}
Fix an object $U$ in $X_{\mr{Zar}}$. To extend $F$, let $V$ be a flat, quasicoherent
$\O_U$-module, and define 
\[
F(V)=\colim_{W\subset V}{F(W)},
\] 
where the colimit is over all coherent $\O_G$-subcomodules
$W\subset V$.  By Lemma \ref{generalized_directed_union}, this is a
direct limit.  Since filtered colimits are exact and commute
with tensor product, $F(V)$ is flat, and the extended functor is 
a tensor functor that is exact.  This establishes (i)

Next, we show that the extended $F$
respects colimits.  Suppose $W=\colim_{\alpha}{W_{\alpha}}$, and write 
$W_{\alpha}=\colim_{\beta}W_{\alpha\beta}$, where each $W_{\alpha\beta}$
is a finite rank $\O_G$-comodule.  Since colimits can be iterated by
\cite{categories}*{IX.8},
we have $W=\colim_{\alpha,\beta}{W_{\alpha,\beta}}$.  It follows that
\[
	F(W)=\colim_{\alpha,\beta}{F(W_{\alpha,\beta})}=
	\colim_{\alpha}\colim_{\beta}F(W_{\alpha\beta})=\colim_{\alpha}F(W_{\alpha}),
\]
hence $F$ respects colimits, which establishes (ii).

It remains to show that $F(\O_G)$ is faithfully
flat.  By \cite{EGAIV}*{2.2.1}, $F(\O_G)$ is faithfully flat over $Y$ if and only if  
the functor $M\mapsto F(\O_G)\otimes_{U'} M$ is an exact and faithful functor
on $\mr{QCoh}_{U'}$ for all $U'\subset Y$ open.  Since $F(\O_G)$ is flat, 
$M\mapsto F(\O_G)\otimes_{U'} M$ is exact.  It remains to show that for any $M\neq 0$,
$F(\O_G)\otimes_{U'} M$ is nonzero.  Since $\O_G$ has $\O_X$ as a direct summand, 
and $F$ is exact, $F(\O_G)$  contains $F(\O_X)=\O_Y$ as a direct summand.  
In particular, $F(\O_G)\otimes_{U'}M=M\oplus M'$ (for some $M'$)
is nonzero, which completes the proof.
\end{proof}

\begin{lemma}\label{extend}
The functor $F$ naturally 
induces a functor from the fibered category over
$X_{\mr{Zar}}$ of $U$-schemes with $G_U$-action that are
flat and affine over $U$ to the fibered category over $X_{\mr{Zar}}$
of schemes flat and affine over
$Y_U$.  The resulting functor, which we again denote by $F$,
respects products and has the property that if $T_0$
has a trivial $G_U$-action then $F(T_0)=Y_U\times_U T_0$.
\end{lemma}

\begin{proof}
Fix an object $U$ in $X_{\mr{Zar}}$.
Let $T$ be a scheme flat and affine over $U$ with $G_U$-action.
Then (the pushforward of) $\O_T$ is an $\O_U$-algebra and
$\O_{G_U}$-comodule.  Furthermore, the multiplication map
$\O_T\otimes\O_T\to\O_T$ is an $\O_{G_U}$-comodule map.
Thus, since $F$ is a tensor functor, 
$F(\O_T)$ is naturally an $\O_{Y_U}$-algebra and flat
by Lemma~\ref{functor_extension}.  We can therefore define
\[
	F(T)=\Spec{F(\O_T)},
\]
a scheme that is flat and affine over $Y_U$. 
Since $F$ is a tensor functor, it is clear that it
respects products.  

To verify the last claim, we identify the full subcategory of trivial
representations in $\op{Rep}_U{G}$ with the category
of finite rank vector bundles on $U$.  For each affine open
$U'\subset Y_U$, we will give a natural isomorphism
\[
	F(V)|U'\iso \O_{U'}\otimes_{\O_U} V
\]
and it will be clear from the construction that these 
isomorphisms will agree on overlaps.  Thus, we may assume
that $Y=\Spec{B}$ is affine.

Furthermore, it suffices to prove the result for 
some affine cover of $X$ so we may assume
that $X=\Spec{A}$ is affine.
Since $F$ is a tensor functor, $F(A)=B$.
Thus, for any vector bundle $V$ the composition
\[
	V\iso\op{Hom}_{A}{(A,V)}\xto{F}\op{Hom}_{B}{(B,F(V))}
\]
gives rise to a natural map of $B$-modules
$\psi: V\otimes_A B\to F(V)$ by adjunction.
Furthermore, $\psi$ is an isomorphism for $V=A^n$.
Since any vector bundle is a direct summand of a free module, 
it follows that $\psi$ is an isomorphism for all $V$.
\end{proof}

\begin{remark}
The above lemma establishes the
aim of this section in the case that $G$ is the trivial group.
When $X$ is not affine, the use of fibered categories is crucial to
establish this result.
\end{remark}

\begin{lemma}\label{F(G)_is_torsor}
Let $P=F(G)$.  Then $P$ is a $G_Y$-torsor naturally in $F$ and $Y$.
\end{lemma}

\begin{proof}
By Lemma~\ref{functor_extension}(iii), $P$ is faithfully flat over $Y$.
Denote the group map $m:G\times G\to G$ and the identity section
$e:X\to G$.  Let $G_0$ denote the same underlying scheme as $G$ with the
trivial $G$-action.
By Lemma~\ref{extend}, applying $F$ to the map of $G$-sets $G\times G_0\to G$ gives rise
to a map $P\times_Y G_Y\to P$.  Again by Lemma~\ref{extend}, 
applying $F$ to the commutative diagrams of $G$-sets

\begin{minipage}{0.4\linewidth}
	\[
		\xymatrix{
			X\times G_0 \ar[r]^{e\times 1} \ar@{=}[dr] & G\times G_0 \ar[d]^{m}\\
			& G_0
		}
	\]
\end{minipage}
\hspace{.5in}
\begin{minipage}{0.4\linewidth}
	\[
		\xymatrix{
			G \times G_0 \times G_0 \ar[r]^-{1\times m} \ar[d]^{m\times 1} 
			& G\times G_0 \ar[d]^{m}\\ 
			G\times G_0 \ar[r]^{m} & G_0
		}
	\]
\end{minipage}
\medskip

\noindent establishes that the map $P\times_Y G_Y\to P$ is a right $G$-action.

Since the $G$-map
\[
	G\times G_0 \to G\times G;\ (g,h)\mapsto (g,gh)
\]
is an isomorphism, the corresponding map induced by $F$,
$P\times_Y G_Y\to P\times_Y P$, is an isomorphism.  
Thus, $P$ is a $G_Y$-torsor. 
\end{proof}

\begin{thm}\label{functors=torsors}
Let $Y$ be faithfully flat scheme over $X$.
The functor from the category of $G_Y$-torsors
to the category of 
tensor functors $F:\op{\mbf{Rep}}{G}\to\mbf{Bun}_Y$ 
that on each fiber over $X_{\mr{Zar}}$ are faithful and exact,
given by
\[
	P\mapsto [F_P: V\mapsto P_U\times^{G_{Y_U}} (V\times_U Y_U)],
\]
is an equivalence of fibered categories.  The quasi-inverse
is given by $F\mapsto F(G)$ (see below remark).
\end{thm}

\begin{remark}\label{F(G)_explanation}
Before we begin the proof, let us summarize the definition of $F(G)$.
As in Lemma \ref{functor_extension}, we can define $F(\O_G)=\colim F(V)$
where $V$ ranges over $\O_X$-coherent $\O_G$-subcomodules of $\O_G$.  Then, as
described in the proof of Lemma \ref{extend}, $F(\O_G)$
is an $\O_X$-algebra, so we can define $F(G)=\Spec{F(\O_G)}$.  The
$G_Y$-action on $F(G)$ is described in the proof
of Lemma \ref{F(G)_is_torsor}.
\end{remark}

\begin{proof}[Proof of Theorem \ref{functors=torsors}]
We must show that the two functors are quasi-inverses.
Given a $G_Y$-torsor $P$, that $F_P(G)$ is naturally
isomorphic to $P$ follows directly from the definition of $F_P$:
\[
	F_P(G) = P\times^G G  = P\times G/[(p,x)\sim(pg,g^{-1}x)] \iso P.
\]
Here the last map is given by $(p,x)\mapsto px$, which respects
the right action on $P\times^G G$ given by $(p,x)\cdot g=(pg,x)=(p,gx)$.

Let $F:\op{\mbf{Rep}}{G}\to\mbf{Bun}_Y$ be given.  Let $P=F(G)$.
We must show that $F_P$ is naturally equivalent to $F$.   
For the remainder of the proof, we will make frequent use of 
Lemma \ref{extend} without explicit mention.  We again use the notation
that if $T$ is some object with $G$-action, then $T_0$ is the same underlying
object with the trivial $G$-action.
Recall that the right $G_Y$-action on $P$ is given by
applying $F$ to the $G$-map $G\times G_0\to G$.  Since $F$ respects products,
$P_U=F(G_U)$ and the right $G_{Y_U}$-action on $P_U$
is given by applying $F$ to $G_U\times_U (G_U)_0\to G_U$.
Fix an object $U$ of $X_{\mr{Zar}}$ and let 
$V$ be a representation of $G_U$.  
Applying $F$ to $\rho: G_U\times_U V_0\to V$ induces a map 
$\phi=F(\rho): P_U\times_{Y_U} (V\times_U Y_U) \to F(V)$.

We first show that $\phi$ factors through the quotient map 
$P_U\times (V\times_U Y_U)\to P_U\times^{G_{Y_U}} (V\times_U Y_U)$.  
By definition, this quotient is defined to be the coequalizer of
\[
	\xymatrix
	{
		P_U\times_{Y_U} G_{Y_U}\times_{Y_U} (V\times_U Y_U)
		\ar@<1ex>[r]^-{\pi_{1,3}} \ar@<-1ex>[r]_-{\beta} &
		P_U\times_{Y_U} (V\times_U Y_U),
	}
\]
where $\beta:(p,g,v)\mapsto (pg,g^{-1}v)$.  Thus, it suffices to show
that $\phi\circ\pi_{1,3}=\phi\circ\beta$.
Denote by $\alpha:G_U\times_U (G_U)_0\times_U V_0\to G_U\times_U V_0$ 
the $G_U$-map $(g,h,v)\mapsto (gh,h^{-1}v)$.
Then it is immediate that the following diagram commutes.
\[
	\begin{CD}
		G_U\times_U (G_U)_0\times_U V_0 @>\pi_{1,3}>> G_U\times_U V_0 \\
		@V{\alpha}VV @VV{\rho}V\\
		G_U\times_U V_0 @>>{\rho}> V
	\end{CD}
\]
By definition of the $G$-action, $\beta=F(\alpha)$.  Thus,  by applying $F$ to
the above diagram, we conclude that $\phi\circ\pi_{1,3}=\phi\circ\beta$.
It follows that $\phi$ descends to a map 
$\phi: P_U\times^{G_{Y_U}} (V\times_U Y_U)\to F(V)$, which 
it remains to show is an isomorphism.

Since $P_U\to Y_U$ is faithfully flat, it suffices to show that $\phi$ is an isomorphism
after pulling back to $P_U$.  One checks from the definitions that we have the 
following sequence of isomorphisms:
\begin{align*}
	P_U\times_{Y_U} (V\times_U Y_U) &\iso (P_U\times_U G_{Y_U})\times^{G_{Y_U}} (V\times_U Y_U)\\
		            &\iso (P_U\times_{Y_U} P_U)\times^{G_{Y_U}} (V\times_U Y_U)\\
		            &\iso P_U\times_{Y_U}(P_U\times^{G_{Y_U}} (V\times_U Y_U)).
\end{align*}
Thus, identifying the source of $1\times\phi$ with the first term in the above
sequence, it remains to show that the induced map 
$\psi:P_U\times_{Y_U} (V\times_U Y_U)\to P_U\times_{Y_U} F(V)$ is an isomorphism.
Following the construction, one sees that $\psi$ comes from applying $F$ 
to the $G_U$-map $G_U\times_U V_0\to G_U\times_U V$ given by $(g,v)\mapsto (g,gv)$.  
Since this latter map is an 
isomorphism, it follows that $\psi$ is an isomorphism, whence the result follows.
\end{proof}

\end{document}